\documentclass[12pt, twoside]{article}
\usepackage{amsmath,amsthm,amssymb}
\usepackage{times}
\usepackage{enumerate}
\usepackage{bm}
\usepackage[all]{xy}
\usepackage{mathrsfs}
\usepackage{amscd}
\usepackage{color}
\def\titlerunning#1{\gdef\titrun{#1}}
\makeatletter
\def\author#1{\gdef\autrun{\def\and{\unskip, }#1}\gdef\@author{#1}}
\def\address#1{{\def\and{\\\hspace*{18pt}}\renewcommand{\thefootnote}{}%
\footnote {#1}}%
\markboth{\autrun}{\titrun}}
\makeatother
\def\email#1{e-mail: #1}
\def\subjclass#1{{\renewcommand{\thefootnote}{}%
\footnote{\emph{Mathematics Subject Classification (2010):} #1}}}
\def\keywords#1{\par\medskip
\noindent\textbf{Keywords.} #1}

\newtheorem{theorem}{Theorem}[section]

\theoremstyle{definition}

\numberwithin{equation}{section}

\input xy
\xyoption{all}

\def\Om{\Omega}

\def\na {\nabla}
\def\Ga{\Gamma}

\begin{document}
\baselineskip=17pt

\titlerunning{A proof on energy gap for Yang-Mills connection}
\title{A proof on energy gap for Yang-Mills connection}

\author{Teng Huang}

\date{}

\maketitle

\address{T. Huang:Key Laboratory of Wu Wen-Tsun Mathematics, Chinese Academy of Sciences,
                   School of Mathematical Sciences, University of Science and Technology of China
                   Hefei, Anhui 230026, PR China; \email{oula143@mail.ustc.edu.cn}}


\begin{abstract}
In this note,\ we prove an ${L^{\frac{n}{2}}}$-energy gap result for Yang-Mills connections on a principal $G$-bundle over a compact manifold without using Lojasiewicz-Simon gradient inequality  (\cite{Feehan2015} Theorem 1.1).
\end{abstract}
\section{Introduction}
Let $X$ be a compact $n$-dimensional Riemannian manifold with smooth Riemannian metric $g$,\ $P\rightarrow X$ a principal $G$-bundle over $X$,\ where $G$ is a compact Lie group.\ We defined the Yang-Mills functional by
$$YM(A)=\int_{X}|F_{A}|^{2}dvol_{g}.$$
where $A$ is a $C^{\infty}$-connection on $P$ and $F_{A}$ is the curvature of $A$.\\
A connection $A$ on $P$ is called Yang-Mills connection,\ if it is a critical point of $YM$ i.e. it obeys the Yang-Mills equation with respect to the metric $g$:
\begin{equation}\label{P2}
d^{\ast}_{A}F_{A}=0
\end{equation}
{In \cite{Feehan2015},\ Feehan proved an $L^{\frac{n}{2}}$-energy gap result for Yang-Mills connections on principal $G$-bundle $P$ over arbitrary closed smooth Riemannian manifold with dimensional $n\geq2$ (\cite{Feehan2015} Theorem 1.1).\ Feehan applied the Lojasiewicz-Simon gradient inequality (\cite{Feehan2015} Theorem 3.2) to remove a positive hypothesis on the Riemannian curvature tensors in a previous $L^{\frac{n}{2}}$-energy gap result due to Gerhardt \cite{Gerhardt} Theorem 1.2.}

{In this note,\ we give another way to prove the $L^{\frac{n}{2}}$-energy gap result of Yang-Mills connection without using the Lojasiewicz-Simon gradient inequality.}
\begin{theorem}\label{T1}({\cite{Feehan2015} Theorem 1.1})
Let $X$ be a compact Riemannian manifold {without boundary} of dimension $n\geq2$ with smooth Riemannan metric $g$,\ $P$ be a $G$-bundle over $X$.\ Then any smooth Yang-Mills connection $A$ over $X$ with compact Lie group $G$ is either satisfies
$${\int_{X}|F_{A}|^{\frac{n}{2}}dvol_{g}\geq C_{0}} $$
for a constant $C_{0}>0$ depending only on $X,n,G$ or the connection $A$ is flat.
\end{theorem}

\section{Preliminaries and basic estimates}
We shall generally adhere to the now standard gauge-theory conventions and notation of Donaldson and Kronheimer \cite{Donaldson/Kronheimer} and Feehan \cite{Feehan2015}.\ Throughout our article,\ $G$ denotes a compact Lie group and $P$ a smooth principal $G$-bundle over a compact Riemannnian manifold $X$ of dimension $n\geq2$ and endowed with Riemannian metric $g$,\ $\mathfrak{g}_{P}$ denote the adjoint bundle of $P$,\ endowed with a $G$-invariant inner product and $\Om^{p}(X,\mathfrak{g}_{P})$ denote the smooth $p$-forms with values in $\mathfrak{g}_{P}$.\ Given a connection on $P$,\ we denote by $\na_{A}$ the corresponding covariant derivative on $\Om^{\ast}(X,\mathfrak{g}_{P})$ induced by $A$ and the Levi-Civita connection of $X$.\ Let $d_{A}$ denote the exterior derivative associated to $\na_{A}$.

For $u\in L^{p}(X,\mathfrak{g}_{P})$,\ where $1\leq p<\infty$ and $k$ is an integer,\ we denote
\begin{equation}\nonumber
\|u\|_{L^{p}_{k,A}(X)}:=\big{(}\sum_{j=0}^{k}\int_{X}|\na^{j}_{A}u|^{p}dvol_{g}\big{)}^{1/p},
\end{equation}
where $\na^{j}_{A}:=\na_{A}\circ\ldots\circ\na_{A}$ (repeated $j$ times for $j\geq0$).\ For $p=\infty$,\ we denote
\begin{equation}\nonumber
\|u\|_{L^{\infty}_{k,A}(X)}:=\sum_{j=0}^{k}ess\sup_{X}|\na^{j}_{A}u|.
\end{equation}
At first,\ {we review a key result due to Uhlenbeck for the connections with $L^{p}$-small curvature ($2p>n$) \cite{Uhlenbeck1985} which provides existence a flat connection $\Ga$ on $P$,\ a global gauge transformation $u$ of $A$ to Coulomb gauge with respect to $\Ga$ and a Sobolev norm estimate for the distance between $\Ga$ and $A$.}
\begin{theorem}({\cite{Uhlenbeck1985} Corollary 4.3 and \cite{Feehan2015} Theorem 5.1})\label{T2}
Let $X$ be a closed,\ smooth manifold of dimension $n\geq2$ and endowed with a Riemannian metric,\ $g$,\ and $G$ be a compact Lie group,\ and $2p>n$.\ Then there are constants,\ $\varepsilon=\varepsilon(n,g,G,p)\in(0,1]$ and $C=C(n,g,G,p)\in[1,\infty)$,\ with the following significance.\ Let $A$ be a $L^{p}_{1}$ connection on a principal $G$-bundle $P$ over $X$.\ If the curvature $F_{A}$ obeying
$$\|F_{A}\|_{L^{p}(X)}\leq\varepsilon,$$
then there exist a flat connection,\ $\Ga$,\ on $P$ and a gauge transformation ${u\in L^{p}_{2}(X)}$ such that\\
(1) $d^{\ast}_{\Ga}\big{(}{u}^{\ast}(A)-\Ga\big{)}=0\ on\ X,$\\
(2) ${\|u^{\ast}(A)-\Ga\|_{L^{p}_{1,\Ga}}\leq C\|F_{A}\|_{L^{p}(X)}}$ and \\
(3) ${\|u^{\ast}(A)-\Ga\|_{L^{\frac{n}{2}}_{1,\Ga}}\leq C\|F_{A}\|_{L^{\frac{n}{2}}(X)}}$.
\end{theorem}
Next,\ {we also review an other key result due to Uhlenbeck concerning an a priori estimate for the curvature of a Yang-Mills connection over a closed Riemannian manifold.}
\begin{theorem}({\cite{Uhlenbeck} Theorem 3.5 and \cite{Feehan2015} Corollary 4.6})\label{T3}
Let $X$ be a compact manifold of dimension $n\geq2$ and endowed with a Riemannian metric $g$, let $A$ be a smooth Yang-Mills connection with respect to the metric $g$ on a smooth $G$-bundle $P$ over $X$.\ Then there exist constants $\varepsilon=\varepsilon(X,n,g)>0$ and $C=C(X,n,g)$ with the following significance.\ If the curvature $F_{A}$ obeying
$${\|F_{A}\|_{L^{\frac{n}{2}}(X)}\leq\varepsilon},$$
then
$$\|F_{A}\|_{L^{\infty}(X)}\leq C\|F_{A}\|_{L^{2}(X)}.$$
\end{theorem}

\section{Proof Theorem \ref{T1}}
For any $p\geq1$ ,\ the estimate in Theorem \ref{T3} yields
\begin{equation}\label{E2.4}
\|F_{A}\|_{L^{p}(X)}\leq C\|F_{A}\|_{L^{\infty}(X)}\leq C\|F_{A}\|_{L^{2}(X)},
\end{equation}
for $C=C(g,n)$.\\
If $n\geq4$,\ by using H\"{o}ler inequality,\ we have
\begin{equation}\label{E2.5}
\|F_{A}\|_{L^{2}(X)}\leq C\|F_{A}\|_{L^{\frac{n}{2}}(X)}.
\end{equation}
If $n=2,3$,\ the $L^{p}$ interpolation implies that
\begin{equation}\nonumber
\begin{split}
\|F_{A}\|_{L^{2}(X)}&\leq C\|F_{A}\|^{\frac{n}{4}}_{L^{\frac{n}{2}}(X)}\|F_{A}\|^{1-\frac{n}{4}}_{L^{\infty}(X)}\\
&\leq  C\|F_{A}\|^{\frac{n}{4}}_{{\frac{n}{2}}(X)}\|F_{A}\|^{1-\frac{n}{4}}_{L^{2}(X)}\\
\end{split}
\end{equation}
and thus
\begin{equation}\label{E2.6}
\|F_{A}\|_{L^{2}(X)}\leq C\|F_{A}\|_{L^{\frac{n}{2}}(X)}.\\
\end{equation}
Therefore,\ by combining (\ref{E2.4})--(\ref{E2.6}),\ we obtain
$$\|F_{A}\|_{L^{p}(X)}\leq C\|F_{A}\|_{L^{\frac{n}{2}}(X)},\ \forall 2p\geq n\ and\ n\geq2.$$
Hence,\ if we suppose $\|F_{A}\|_{L^{\frac{n}{2}}(X)}$ sufficiently small such that $\|F_{A}\|_{L^{q}(X)}$ ($2q>n$ and $n\geq2$) satisfies the hypothesis of Theorem \ref{T2},\ then Theorem \ref{T2} provides a flat connection $\Ga$ on $P$,\ and a gauge transformation $u\in\mathcal{G}_{P}$ and the estimate
$$\|u^{\ast}(A)-\Ga\|_{L^{q}_{1}(X)}\leq C(q)\|F_{A}\|_{L^{q}(X)},$$
and
$$d_{\Ga}^{\ast}(u^{\ast}(A)-\Ga)=0.$$
We denote ${\tilde{A}:=u^{\ast}(A)}$ and  ${a:=u^{\ast}(A)-\Ga}$,\ then the curvature of $\tilde{A}$ is
$$F_{\tilde{A}}=d_{\Ga}a+a\wedge a.$$
The connection $\tilde{A}$ also satisfies Yang-Mills equation
\begin{equation}\label{G4}
0=d_{\tilde{A}}^{\ast}F_{\tilde{A}}.
\end{equation}
Hence taking the $L^{2}$-inner product of (\ref{G4}) with $a$,\ we obtain
\begin{equation}\nonumber
\begin{split}
0&=(d_{\tilde{A}}^{\ast}F_{\tilde{A}},a)_{L^{2}(X)}\\
&=(F_{\tilde{A}},d_{\tilde{A}}a)_{L^{2}(X)}\\
&=(F_{\tilde{A}},d_{\Ga}a+2a\wedge a)_{L^{2}(X)}\\
&=(F_{\tilde{A}},F_{\tilde{A}}+a\wedge a)_{L^{2}(X)}.\\
\end{split}
\end{equation}
Then we get
\begin{equation}\nonumber
\begin{split}
\|F_{A}\|^{2}_{L^{2}(X)}&=\|F_{\tilde{A}}\|^{2}_{L^{2}(X)}\\
&=-(F_{\tilde{A}},a\wedge a)_{L^{2}(X)}\\
&\leq\|F_{\tilde{A}}\|_{L^{2}(X)}\|a\wedge a\|_{L^{2}(X)}\\
&=\|F_{A}\|_{L^{2}(X)}\|a\wedge a\|_{L^{2}(X)}\\
\end{split}
\end{equation}
here we use the fact {$|F_{u^{\ast}(A)}|=|F_{A}|$ since $F_{u^{\ast}(A)}=u\circ F_{A}\circ u^{-1}$}.\\
If $n\geq4$:\\
\begin{equation}\nonumber
\begin{split}
\|a\wedge a\|_{L^{2}(X)}&\leq C\|a\|^{2}_{L^{4}(X)}\\
&\leq C\|a\|^{2}_{L^{n}(X)}\\
&\leq C\|a\|^{2}_{L^{\frac{n}{2}}_{1}(X)}\\
&{\leq C\|F_{A}\|^{2}_{L^{\frac{n}{2}}(X)}}\\
&\leq C\|F_{A}\|^{2}_{L^{\infty}(X)}\\
&\leq C\|F_{A}\|^{2}_{L^{2}(X)}.\\
\end{split}
\end{equation}
where we apply the Sobolev embedding $L^{\frac{n}{2}}_{1}\hookrightarrow L^{n}$.\\
If $n=2,3$,\
\begin{equation}\nonumber
\begin{split}
\|a\wedge a\|_{L^{2}(X)}&\leq C\|a\|^{2}_{L^{4}(X)}\\
&{\leq C\|a\|^{2}_{L^{2}_{1}(X)}}\\
&\leq C\|F_{A}\|^{2}_{L^{2}(X)}.\\
\end{split}
\end{equation}
{where we apply the Sobolev embedding $L^{2}_{1}\hookrightarrow L^{4}$.}\\
By combining the preceding inequalities we have
$$\|F_{A}\|^{2}_{L^{2}(X)}\leq C\|F_{A}\|^{3}_{L^{2}(X)},$$
We can choose $\|F_{A}\|_{L^{2}(X)}$ sufficiently small to such that $C\|F_{A}\|_{L^{2}(X)}<1$,\ hence {$\|F_{A}\|_{L^{2}(X)}\equiv0$ and thus $A$ must be a flat connection}.\ Then we complete the proof.\\
\subsection*{Acknowledgements}
I would like to thank Professor Paul Feehan for helpful comments regarding his article \cite{Feehan2015}.\ { I thank the anonymous
referee for a careful reading of my article and helpful comments and corrections.}\ This work is partially supported by Wu Wen-Tsun Key Laboratory of Mathematics of Chinese Academy of Sciences at USTC.


\bigskip
\footnotesize
\begin{thebibliography}{SK}


\bibitem{Donaldson/Kronheimer}
Donaldson S.~K., Kronheimer P.~B.,
The geometry of four-manifolds,
\textit{Oxford University Press}, 1990.
\bibitem{Feehan2015}
Feehan P.~M.~N.,
Energy gap for Yang-Mills connections, II: Arbitrary closed Riemannian manifolds,
\textit{Adv. Math.} Doi.org/10.1016/j.aim.2017.03.023
\bibitem{Gerhardt}
Gerhardt C.,
An energy gap for Yang-Mills connections,
\textit{Comm. Math. Phys.} \textbf{298} (2010), 515--522.

\bibitem{Uhlenbeck}
Uhlenbeck K.~K.,
Removable singularites in Yang-Mills fields,
\textit{Comm. Math. Phys.} \textbf{83} (1982), 11--29.

\bibitem{Uhlenbeck1985}
Uhlenbeck K.~K.,
The Chern classes of Sobolev connections,
\textit{Comm. Math. Phys.} \textbf{101} (1985), 445--457.


\baselineskip=17pt
\end{thebibliography}
\end{document}